\documentclass[10pt]{article}
\usepackage{amsmath}
\usepackage{amsfonts}
\usepackage{amssymb}
\usepackage{amsthm}

\begin{document}

\newcommand{\C}{{\mathbb{C}}}
\newcommand{\R}{{\mathbb{R}}}
\newcommand{\Z}{{\mathbb{Z}}}
\newcommand{\q}{\left}
\newcommand{\w}{\right}
\newcommand{\lap}{\bigtriangleup}
\newcommand{\grad}{\bigtriangledown}
\newcommand{\chO}{\mathrm{ch}\q(\Omega\w)}
\newcommand{\sE}{\mathcal{E}}
\newcommand{\sPqnp}{\mathcal{P}_q}
\newcommand{\sPznp}{\mathcal{P}_0}
\newcommand{\sPq}[1]{\sPqnp\q( #1\w)}
\newcommand{\sPz}[1]{\sPznp\q( #1\w)}
\newcommand{\LtO}{L^{2}\q( \Omega \w)}
\newcommand{\HoO}{H^1\q( \Omega\w)}
\newcommand{\HtO}{H^2\q( \Omega\w)}
\newcommand{\HozO}{H^1_0\q( \Omega\w)}
\newcommand{\HohdO}{H^{\frac{1}{2}}\q( \partial \Omega\w)}
\newcommand{\HmohdO}{H^{-\frac{1}{2}}\q( \partial \Omega\w)}
\newcommand{\HmthdO}{H^{-\frac{3}{2}}\q( \partial \Omega\w)}
\newcommand{\HthdO}{H^{\frac{3}{2}}\q( \partial \Omega\w)}
\newcommand{\HDO}{H_{\lap}\q( \Omega\w)}
\newcommand{\Ft}{\widetilde{F}}
\newcommand{\Bt}{\widetilde{B}}
\newcommand{\Btt}{\widetilde{\widetilde{B}}}
\newcommand{\Gc}{\check{\Gamma}}
\newcommand{\zb}{\overline{z}}
\newcommand{\bq}{b_q}
\newcommand{\bz}{b_{0}}
\newcommand{\sH}{\mathcal{H}\q( \partial \Omega\w)}
\newcommand{\tr}[1]{{\mathrm{tr}}\q( #1\w)}
\newcommand{\trnp}{{\mathrm{tr}}}

\newcommand{\LtOn}[1]{\q\| #1 \w\|_{\LtO}}

\newcommand{\Opp}[3]{{#1}_{#3}}
\newcommand{\Opt}[2]{\Opp{#1}{#2}{\tau}}
\newcommand{\Opmt}[2]{\Opp{#1}{#2}{-\tau}}
\newcommand{\Optb}[2]{\Opp{\overline{#1}}{#2}{\tau}}
\newcommand{\Opmtb}[2]{\Opp{\overline{#1}}{#2}{-\tau}}
\newcommand{\sLt}{\Opt{\mathcal{L}}{W}}
\newcommand{\sLtb}{\Optb{\mathcal{L}}{W}}
\newcommand{\sLmt}{\Opmt{\mathcal{L}}{W}}
\newcommand{\sLmtb}{\Opmtb{\mathcal{L}}{W}}
\newcommand{\Gt}{\Opt{G}{W}}
\newcommand{\Gmt}{\Opmt{G}{W}}
\newcommand{\Gtb}{\Optb{G}{W}}
\newcommand{\Gmtb}{\Opmtb{G}{W}}
\newcommand{\Ht}{\Opt{H}{W}}
\newcommand{\Hmt}{\Opmt{H}{W}}
\newcommand{\Htb}{\Optb{H}{W}}
\newcommand{\Hmtb}{\Opmtb{H}{W}}
\newcommand{\Tt}{\Opt{T}{W}}
\newcommand{\Tmt}{\Opmt{T}{W}}
\newcommand{\Ttb}{\Optb{T}{W}}
\newcommand{\Tmtb}{\Opmtb{T}{W}}
\newcommand{\St}{S_\tau}

\renewcommand{\Im}{\mathrm{Im}}
\newcommand{\LtdOp}{L^2\q( \partial \Omega_{+} \w)}
\newcommand{\LtdOm}{L^2\q( \partial \Omega_{-} \w)}
\newcommand{\LtdOpn}[1]{\q\| #1 \w\|_{\LtdOp}}
\newcommand{\LtdOmn}[1]{\q\| #1 \w\|_{\LtdOm}}
\newcommand{\Ltnp}[2]{\q\| #1 \w\|_{L^2\q( #2 \w)}}
\newcommand{\sgn}[1]{{\mathrm{sgn}}\q( #1 \w)}
\newcommand{\LtOOpn}[1]{\q\| #1 \w\|_{\LtO\rightarrow \LtO}}
\newcommand{\sD}{\mathcal{D}}
\newcommand{\sT}{\mathcal{T}}
\newcommand{\pt}{\pi_{\tau}}
\newcommand{\ptb}{\overline{\pi}_{\tau}}
\newcommand{\pmt}{\pi_{-\tau}}
\newcommand{\pmtb}{\overline{\pi}_{-\tau}}
\newcommand{\ptq}{\pi_{\tau}^q}
\newcommand{\ptbq}{\overline{\pi}_{\tau}^q}
\newcommand{\pmtq}{\pi_{-\tau}^q}
\newcommand{\pmtbq}{\overline{\pi}_{-\tau}^q}
\newcommand{\ptz}{\pi_{\tau}^0}
\newcommand{\ptbz}{\overline{\pi}_{\tau}^0}
\newcommand{\pmtz}{\pi_{-\tau}^0}
\newcommand{\pmtbz}{\overline{\pi}_{-\tau}^0}
\newcommand{\LtOip}[2]{\q< #1,#2 \w>_{\LtO}}
\newcommand{\LtdOip}[2]{\q< #1,#2 \w>_{\LtdO}}
\newcommand{\LtdOmip}[2]{\q< #1,#2 \w>_{\LtdOm}}
\newcommand{\LtdOpip}[2]{\q< #1,#2 \w>_{\LtdOp}}
\newcommand{\Ltipp}[3]{\q< #1,#2 \w>_{#3}}
\newcommand{\Mt}{M_{\tau}}
\newcommand{\sMt}{\mathcal{M}_\tau}
\newcommand{\Rt}{R_\tau}
\newcommand{\Rmt}{R_{-\tau}}

\newcommand{\mut}{\mu_\tau}
\newcommand{\nut}{\nu_\tau}
\newcommand{\ut}{u_\tau}
\newcommand{\vt}{v_\tau}
\newcommand{\numt}{\nu_{-\tau}}
\newcommand{\vmt}{v_{-\tau}}
\newcommand{\wt}{w_\tau}
\newcommand{\ot}{\omega_\tau}
\newcommand{\htt}{\widetilde{h}}

\newtheorem{thm}{Theorem}[section]
\newtheorem{cor}[thm]{Corollary}
\newtheorem{prop}[thm]{Proposition}
\newtheorem{lemma}[thm]{Lemma}
\newtheorem{conj}[thm]{Conjecture}


\theoremstyle{remark}
\newtheorem{rmk}[thm]{Remark}

\theoremstyle{definition}
\newtheorem{defn}[thm]{Definition}

\title{Reconstruction in the Calder\'on Problem with Partial Data}
\author{Adrian Nachman\footnote{nachman@math.toronto.edu; Department of Mathematics, University of Toronto, Room 6290, 40 St. George Street, Toronto, Ontario; The Edward S. Rogers Sr. Department of Electrical and Computer Engineering, University of Toronto, 10 King's College Road, Toronto, Ontario}{} {} and Brian Street\footnote{street@math.wisc.edu; University of Wisconsin-Madison, Department of Mathematics, 480 Lincoln Dr., 53706}}
\date{}

\maketitle

\begin{abstract}
We consider the problem of recovering the coefficient $\sigma\q( x\w)$ of
the elliptic equation $\grad \cdot\q(\sigma \grad u\w)=0$ in a body
from measurements of the Cauchy data on possibly very small
subsets of its surface.  We give a constructive proof of a uniqueness
result by Kenig, Sj\"ostrand, and Uhlmann.  We construct a
uniquely specified family of solutions such that their traces
on the boundary can be calculated by solving an integral
equation which involves only the given partial Cauchy data.
The construction entails a new family of Green's functions
for the Laplacian, and corresponding single layer potentials,
which may be of independent interest.
\end{abstract}

\section{Introduction}
Let $\Omega$ be a bounded domain in $\R^n$, $n\geq 3$, with $C^2$ boundary,
and let $\sigma$ be a strictly positive function in $C^2\q( \bar{\Omega}\w)$.
The Dirichlet-to-Neumann map is the operator on the boundary
$\Lambda_\sigma:\HohdO\rightarrow \HmohdO$ defined as
$$\Lambda_\sigma f = \sigma \partial_\nu u\big|_{\partial\Omega},$$
where $u\in \HoO$ is the solution of the Dirichlet problem:
\begin{equation}\label{EqnDirichletProblem}
\grad\cdot\q( \sigma \grad u\w)=0\text{ in }\Omega, \quad u\big|_{\partial\Omega} =f
\end{equation}
and $\nu$ denotes the exterior unit normal to $\Omega$.
If $\Omega$ models an inhomogeneous, isotropic body with conductivity $\sigma$
then $\Lambda_\sigma f$ is the normal component of the current flux at the boundary 
corresponding to a voltage potential $f$ on $\partial \Omega$.

In 1980, Calder\'on \cite{CalderonOnAnInverseBoundaryValueProblem} posed
the following problem:  decide whether $\sigma$ is uniquely determined
by $\Lambda_\sigma$ and, if so, find a method to reconstruct
$\sigma$ from knowledge of $\Lambda_\sigma$.  The problem is of
practical interest in medical imaging and geophysics, where one
seeks to image the conductivity of a body
by making voltage and current measurements
at its surface.  For a summary of the considerable
progress achieved on Calder\'on's problem since his groundbreaking paper,
see \cite{GKLU}, Section 2.

Recent work has shown that uniqueness in the above problem holds even 
if measurements are available only on part of the boundary.
Bukhgeim and Uhlmann \cite{BukhgeimUhlmannRecoveringAPotentialFromPartialCauchyData}
proved that knowledge of values of $\Lambda_\sigma$ on, roughly, slightly
more that half of the boundary $\partial\Omega$ for all $f$ uniquely
determines the conductivity $\sigma$ in $\Omega$ (assuming it is known
on $\partial\Omega$).  This was improved by Kenig, Sj\"ostrand,
and Uhlmann \cite{KenigSjostrandUhlmannTheCalderonProblemWithPartialData}
who assumed $\Lambda_\sigma f$ known on a possibly very small open
subset $U$ of the boundary for $f$
supported in a neighborhood of $\partial\Omega\setminus U$.  (We
describe this result more precisely below.)

The methods in \cite{KenigSjostrandUhlmannTheCalderonProblemWithPartialData}
are non-constructive:  one assumes that one is given two
Dirichlet-to-Neumann maps which agree on appropriate
subsets of the boundary and one shows that the corresponding conductivities
must also agree.  In this paper, we give a reconstruction method.
As in the solution of the the reconstruction
part of Calder\'on's problem in \cite{N1}, we would like to set up
an integral equation on the boundary which in this case
involves only the given data and yields the boundary
values of the geometric optics solutions introduced in \cite{KenigSjostrandUhlmannTheCalderonProblemWithPartialData}.
The main difficulty is that the complex geometrical optics solutions
of \cite{KenigSjostrandUhlmannTheCalderonProblemWithPartialData}
are highly non-unique.  Starting from the Carleman estimate
of \cite{KenigSjostrandUhlmannTheCalderonProblemWithPartialData}
we show how to construct new solutions which are uniquely
specified and for which the boundary values
can be calculated by solving an integral equation which involves only
the assumed partial knowledge of the Cauchy data.
To do so we construct, given a (possibly small) open
subset $U$ of $\partial \Omega$ as above, a new family
of Green's functions $G\q( x,y\w)$ for the Laplacian which vanish, roughly
speaking, when $x\in U$ or when $y\in \partial\Omega\setminus U$ (see
Theorem \ref{ThmGreensThm} for a precise statement).
We also give a novel treatment of the boundedness properties of the corresponding
single layer operators, which may be of independent interest.
These are the main ingredients needed for our boundary
integral equation.

We now turn to more rigorous details.
Fix any point $x_0$ in $\R^n\setminus \overline{\chO}$,
the complement of the closure of the convex hull, $\chO$, of
$\Omega$.  Following \cite{KenigSjostrandUhlmannTheCalderonProblemWithPartialData},
we define the front and back faces of $\partial\Omega$ by
\begin{equation*}
\begin{split}
F\q( x_0\w) &= \q\{x\in \partial\Omega: \q( x-x_0\w)\cdot \nu\q( x\w) \leq 0\w\}\\
B\q( x_0\w) &= \q\{x\in \partial\Omega: \q( x-x_0\w)\cdot \nu\q( x\w) \geq 0\w\}.
\end{split}
\end{equation*}
The uniqueness result of \cite{KenigSjostrandUhlmannTheCalderonProblemWithPartialData}
can then be stated as follows:
\begin{thm}[\cite{KenigSjostrandUhlmannTheCalderonProblemWithPartialData}, Cor. 1.4]
Let $\Omega$, $x_0$, $F\q( x_0\w)$, and $B\q( x_0\w)$ be as above,
and let $\sigma_1, \sigma_2\in C^2\q( \overline{\Omega}\w)$ be strictly
positive.  Assume that $\sigma_1=\sigma_2$ on $\partial\Omega$.
Suppose that there exist open neighborhoods $\widetilde{F}, \widetilde{B}\subset \partial\Omega$
of $F\q( x_0\w)$ and $B\q( x_0\w)$ respectively, such that $\Lambda_{\sigma_1}f = \Lambda_{\sigma_2}f$
in $\widetilde{F}$ for all $f\in \HohdO$ supported in $\widetilde{B}$.  Then
$\sigma_1=\sigma_2$ in $\Omega$.
\end{thm}

The above theorem was obtained in \cite{KenigSjostrandUhlmannTheCalderonProblemWithPartialData}
as a consequence of the following result for Schr\"odinger operators.
Let $q\in L^\infty\q( \Omega\w)$ (possibly complex valued), and
assume that $0$ is not a Dirichlet eigenvalue of $-\lap +q$ in $\Omega$.
Then for any $v\in \HohdO$ there is a unique (weak) solution $w\in \HoO$
of
\begin{equation}\label{EqnDirichletProblemSch}
\q( -\lap +q\w) w=0 \text{ in } \Omega
\end{equation}
with $w\big|_{\partial_\Omega}=v$.  Define the corresponding Dirichlet-to-Neumann
map $\Lambda_q: \HohdO\rightarrow \HmohdO$ by
\begin{equation}\label{EqnDTNDefn}
\q<\tr{w_0}, \Lambda_q v\w> = \int_\Omega \grad w_0 \cdot \grad w + q w_0 w \text { for any } w_0\in \HoO,
\end{equation}
where $\q<\cdot,\cdot\w>$ denotes the bilinear paring of $\HohdO$ and $\HmohdO$.

\begin{thm}[\cite{KenigSjostrandUhlmannTheCalderonProblemWithPartialData}, Theorem 1.1]\label{ThmKSU}
Let $x_0$, $\Omega$, $F\q( x_0\w)$, and $B\q( x_0\w)$ be as above,
and let $q_i\in L^\infty\q( \Omega\w)$, $i=1,2$, be two potentials such that
$0$ is not a Dirichlet eigenvalue of $-\lap+q_i$ in $\Omega$.  Suppose
that there exist open neighborhoods $\widetilde{F},\widetilde{B}\subset \partial\Omega$
of $F\q( x_0\w)$ and $B\q( x_0\w)$ respectively, such that
$\Lambda_{q_1} v=\Lambda_{q_2} v$ in $\widetilde{F}$,
for all $v\in \HohdO$ with support in $\widetilde{B}$.
Then, $q_1=q_2$.  
\end{thm}

The well-known substitution $u=\sigma^{-1/2} w$
in \eqref{EqnDirichletProblem}
yields a solution $w$ of \eqref{EqnDirichletProblemSch}, with
$q=\frac{\lap \sigma^{1/2}}{\sigma^{1/2}}$, and
\begin{equation}\label{EqnSchLambdaFromCon}
\Lambda_q = \sigma^{-1/2} \q(\Lambda_\sigma + \frac{1}{2} \frac{\partial\sigma}{\partial\nu}\big|_{\partial\Omega}\w)\sigma^{-1/2}.
\end{equation}

We note that, for $\sigma\in C^2\q( \overline{\Omega}\w)$, with $\sigma\big|_{\partial\Omega}$ known, $\frac{\partial\sigma}{\partial\nu}$ can
be reconstructed on $\widetilde{F}\cap \widetilde{B}$ from
measurements of $\Lambda_\sigma f$ on $\widetilde{F}$ for all
$f\in \HohdO$ supported in $\widetilde{B}$.  (See Theorem 6(ii) in \cite{N2}.)
Thus, we henceforth assume known the map
$v\mapsto \Lambda_q v\big|_{\widetilde{F}}$
for $v$ supported in $\Bt$ (see Remark \ref{RmkFuncSpace} for the precise
class of $v$).

As mentioned earlier, the proof of Theorem \ref{ThmKSU} in \cite{KenigSjostrandUhlmannTheCalderonProblemWithPartialData}
is nonconstructive, and begins with the assumption that one is given two
such $q_1$ and $q_2$ for which the partial boundary data agree.
Under these assumptions, it was shown in \cite{DosSantosFerreiraKenigSjostrandUhlmannDeterminingAMagneticSchrodingerOperatorFromPartialCauchyData}
that one can conclude that certain Radon transform information
of $q_1-q_2$ must vanish, and this is enough to
show that $q_1=q_2$ (actually, \cite{DosSantosFerreiraKenigSjostrandUhlmannDeterminingAMagneticSchrodingerOperatorFromPartialCauchyData}
deals with more general magnetic Schr\"odinger operators).
The goal of this paper is to show how, given the map
$v\mapsto \Lambda_q v\big|_{\widetilde{F}}$ for $v$ supported
in $\Bt$,
one may reconstruct the aforementioned Radon transform information of $q$.

We now describe more precisely the transform our method reconstructs.
We follow the presentation of \cite{DosSantosFerreiraKenigSjostrandUhlmannDeterminingAMagneticSchrodingerOperatorFromPartialCauchyData}, which provides a
change of variables which will simplify the exposition.  Fix $R>0$ so large that
$\overline{\Omega}\subset B\q( x_0,R\w)$, let $H$ be a hyperplane separating
$x_0$ and $\chO$, and let $H^{+}$ denote the corresponding open half space containing
$\overline{\Omega}$.  Set
$$\Gamma = \q\{ \theta\in S^{n-1} : x_0+R\theta\in H^{+} \w\}$$
and let $\Gc$ denote the image of $\Gamma$ under the antipodal map.
Fix $\alpha_0\in S^{n-1}\setminus \q( \Gamma \cup \Gc\w)$.
It is important that both $x_0$ and $\alpha_0$ may be perturbed
slightly, and all of our assumptions remain intact.

With this $x_0$ and $\alpha_0$ fixed, we may translate and rotate
$\Omega$ so that, without loss of generality, $x_0=0$ and
$\alpha_0=\q( 1,0,\ldots,0\w)$; note, then,
that $\overline{\Omega}$ does not intersect the line $\R\times \q\{ 0\w\}\times\cdots \times \q\{ 0\w\} \subset \R^n$.

For $x\in \R^n$, we write $x=\q( x_1, x'\w) \in \R\times \R^{n-1}$.  We
then switch to polar coordinates in the $x'$ variable.  Indeed,
denote by $\q( x_1, r, \theta\w)\in \R\times \R_{+}\times S^{n-2}$
such a coordinate system.  Note that, since $\overline{\Omega}$ does
not intersect $\R\times \q\{0\w\}\times \cdots \times \q\{0\w\}$, these
coordinates are good on all of $\overline{\Omega}$.  Let $z$ denote
the complex variable $z=x_1+ir$.  We have:
\begin{equation}\label{EqnLapEqn}
\begin{split}
\lap &= \frac{\partial^2}{\partial x_1^2} + \frac{\partial^2}{\partial r^2} + \frac{n-2}{r} \frac{\partial}{\partial r} + \frac{1}{r^2} \lap_{S^{n-2}}\\
&= 4 \frac{\partial^2}{\partial z \partial\zb} +\frac{2\q( n-2\w)}{z-\zb}\q(\frac{\partial}{\partial\zb}-\frac{\partial}{\partial z}\w) + \frac{1}{\q(z-\zb\w)^2}\lap_{S^{n-2}}.
\end{split}  
\end{equation}

With this notation, we now state our main result.
\begin{thm}\label{ThmMainThm}
Let $x_0$, $\Omega$, $F\q( x_0\w)$, $B\q( x_0\w)$ be as above, and let
$\Ft, \Bt\subset \partial\Omega$ be open neighborhoods of $F\q( x_0\w)$,
respectively $B\q( x_0\w)$.  Let $q\in L^\infty\q( \Omega\w)$ be such that
$0$ is not a Dirichlet eigenvalue of $-\lap+q$ in $\Omega$.  Given
$\Lambda_q v$ on $\Ft$ for all\footnote{See Remark \ref{RmkFuncSpace} for the precise class of $v$ we work with.} $v$ supported in $\Bt$ one can
reconstruct the integrals
\begin{equation}\label{EqnIntToRecon}
\int q\q( x_1, r, \theta\w) g\q( \theta\w)\: dx_1\: dr\: d\theta
\end{equation}
for all $g\in C^\infty\q( S^{n-2}\w)$.  (Here $d\theta$ denotes the usual surface measure on the unit sphere $S^{n-2}$.)
\end{thm}

By varying $x_0$ and $\alpha_0$ slightly (staying within the given data),
it is shown in \cite{DosSantosFerreiraKenigSjostrandUhlmannDeterminingAMagneticSchrodingerOperatorFromPartialCauchyData}
that the resulting integrals determine $q$; we refer the reader to that
paper for the details of the proof.

A brief outline of our paper is as follows.  In Section \ref{SectionFunctionSpaces}
we define the function space on the boundary in which our integral equation
will be solved.  In Section \ref{SectionGreensOp} we construct
the new Green's operators for the Laplacian.  In Section \ref{SectionSpecialSolutions}
we select appropriate uniquely specified complex geometrical optics solutions from those of \cite{KenigSjostrandUhlmannTheCalderonProblemWithPartialData}
when $q=0$.  These will serve as ``incident waves'' in our construction.
In Section \ref{SectionSpecialSolutionsq} we define our new solutions
and the corresponding nonlinear transform $t\q( \tau, q\w)$ of $q$.
In Section \ref{SectionIntegralEqns}, we introduce the new single layer
operators and prove the unique solvability of our boundary integral equation.
This yields the reconstruction of $t\q( \tau, q\w)$ from the partial data,
and the proof of Theorem \ref{ThmMainThm}.

\section{Function Spaces}\label{SectionFunctionSpaces}
We define the Bergman space
$$\bq=\q\{u\in \LtO: \q(-\lap+q\w)u=0\w\}$$
and topologize it as a closed subspace of $\LtO$.  We define
the harmonic Bergman space $\bz$ in a similar way, with $q$ replaced
by $0$.

Following \cite{BukhgeimUhlmannRecoveringAPotentialFromPartialCauchyData},
we work with the Hilbert space
$$\HDO = \q\{u\in \LtO: \lap u\in \LtO\w\},$$
the maximal domain of the Laplacian, with norm
$$\q\|u\w\|_{\HDO}^2 = \LtOn{u}^2+\LtOn{\lap u}^2.$$
The trace map
$$\tr{u} = u\big|_{\partial\Omega}$$
extends to a continuous map $\HDO\rightarrow \HmohdO$; moreover,
if $u\in \HDO$ and $\tr{u}\in \HthdO$, then
$u\in \HtO$ (see \cite{BukhgeimUhlmannRecoveringAPotentialFromPartialCauchyData}, \cite{LM}).
We define:
$$\sH = \q\{\tr{u} : u\in \HDO\w\}\subset \HmohdO$$
though for the moment, we do not define a topology on $\sH$.
The space $\sH$ will be the setting for our main boundary integral equation.

Note that $\bq\subset \HDO$, and so the trace map makes sense as a map
$\bq\rightarrow \sH$.  In fact, this map is
one-to-one and onto.
\begin{prop}\label{PropTraceMaps}
If $q\in L^\infty\q( \Omega\w)$ and $0$ is not a Dirichlet eigenvalue of
$-\lap+q$ in $\Omega$, then the trace map $\trnp:\bq\rightarrow \sH$
is one-to-one and onto.
\end{prop}
\begin{proof}
Suppose $u,v\in \bq$, with $\tr{u}=\tr{v}$.  Then $w=u-v\in \bq$,
with $\tr{w}=0$.  Hence, $w\in \HtO$, $\q( -\lap+q\w)w=0$, so by
the hypothesis on $q$, $w=0$.  Thus $\trnp$ is one-to-one.

Suppose $g\in \sH$.  Thus there exists a function $u\in \HDO$
such that $\tr{u}=g$.  Let $v$ be the $\HozO$ solution to the
Dirichlet problem $\q( -\lap+q\w)v=\q( -\lap+q\w) u$, $\tr{v}=0$,
and let $w=u-v$.  Then $w\in \HDO$, $\tr{w}=g$ and $\q( -\lap+q\w) w=0$.
Thus $\trnp$ is onto.
\end{proof}

We define $\sPznp$ to the the inverse of $\trnp:\bz\rightarrow \sH$ and
$\sPqnp$ to be the inverse of $\trnp: \bq\rightarrow \sH$.
We now define the norm on $\sH$ by $\q\|g\w\|_{\sH}=\LtOn{\sPz{g}}$.
With this topology on $\sH$, the above maps are all continuous.
\begin{lemma}
The map $\trnp: \HDO\rightarrow \sH$ is continuous, and under the
hypothesis of Proposition \ref{PropTraceMaps}, $\trnp:\bq\rightarrow \sH$
is a homeomorphism.
\end{lemma}
\begin{proof}
Take $u\in \HDO$, let $v$ be the unique $H^{1}_0\q( \Omega\w)$ solution
to the Dirichlet problem $\lap v=\lap u$, $\tr{v}=0$, and let $w=u-v$.
Note that the map $u\mapsto w$ is continuous $\HDO\rightarrow \LtO$,
and since $\lap w=0$, $u\mapsto w$ is continuous $\HDO \rightarrow \bz$.
Thus, $u\mapsto \tr{w}$ is continuous $\HDO\rightarrow \sH$; however
$\tr{u}=\tr{w}$, establishing the first claim.

Since $\bq$ continuously embeds into $\HDO$, we have that
$\trnp:\bq\rightarrow \sH$ is continuous.  Since it is bijective
(Proposition \ref{PropTraceMaps}) the open mapping theorem shows
that it is a homeomorphism.
\end{proof}

Having extended the solvability of the Dirichlet problem
to boundary data in $\sH$, we now turn to the Dirichlet-to-Neumann map.

\begin{prop}\label{PropExtendDTN}
Assume $q\in L^\infty\q( \Omega\w)$ and $0$ is not a Dirichlet
eigenvalue of $-\lap+q$ in $\Omega$.  Then $\Lambda_q-\Lambda_0$
extends to a continuous map $\sH\rightarrow \sH^{*}$.
\end{prop}
\begin{proof}
Suppose $f,g\in \HohdO$.  From \eqref{EqnDTNDefn} we have
\begin{equation}\label{EqnAlles}
\q<g, \q(\Lambda_q-\Lambda_0 \w)f\w>=\int_\Omega \sPz{g} q\sPq{f}.
\end{equation}
The right hand side extends continuously to all $f,g\in \sH$ and therefore
so does the left hand side.
\end{proof}

\begin{rmk}\label{RmkFuncSpace}
We will henceforth assume knowledge of $\q( \Lambda_q-\Lambda_0\w)f\big|_{\Ft}$
for $f\in \sH\cap \mathcal{E}'\q(\Bt\w)$.
\end{rmk}

\section{The Green's Operators}\label{SectionGreensOp}
For $\tau\in \R$ (later on we will take $\q|\tau\w|$ large), define:
\begin{equation*}
\sLt = z^{-\tau} \lap z^{\tau}, \quad \sLtb = \zb^{-\tau} \lap \zb^{\tau}.
\end{equation*}
Here, $\tau$ is playing the role that $\frac{1}{h}$ played in
\cite{KenigSjostrandUhlmannTheCalderonProblemWithPartialData,DosSantosFerreiraKenigSjostrandUhlmannDeterminingAMagneticSchrodingerOperatorFromPartialCauchyData}.
Note that since $\overline{\Omega}$ lies in the open half plane $r=\Im \q( z\w) >0$,
$z^{\tau}\in C^{\infty}\q( \overline{\Omega}\w)$ (where $z^{\tau}$ is defined
via the principal branch of the logarithm).

\begin{rmk}
Of course, $\sLmtb$ is the formal adjoint of $\sLt$, however we use the above notation since
we will construct Green's operators for $\sLt$ and $\sLmtb$ in tandem,
which we will not {\it a priori} know to be adjoints of each other.
\end{rmk}

We define:
\begin{equation*}
\partial \Omega_{\pm} = \q\{ x\in \partial \Omega : \pm x\cdot \nu\q( x\w) \geq 0 \w\}
\end{equation*}
Note that $F\q( 0 \w) = \partial\Omega_{-}$ and $B\q( 0\w) = \partial\Omega_{+}$.
For $x\in \partial \Omega$, we define $\gamma\q( x\w) = \frac{\sqrt{\q| x\cdot \nu\q( x\w) \w|}}{\q| x\w|}$.
Then, the
Carleman estimate of \cite{KenigSjostrandUhlmannTheCalderonProblemWithPartialData} (see also \cite{DosSantosFerreiraKenigSjostrandUhlmannDeterminingAMagneticSchrodingerOperatorFromPartialCauchyData})
can be written as: for $\q|\tau\w|>0$ and all $u\in C^{\infty}\q( \overline{\Omega}\w)$ with $\tr{u}=0$,
\begin{equation}\label{EqnCarleman}
\begin{split}
&\q|\tau\w|^{-\frac{1}{2}}\Ltnp{\gamma \partial_{\nu} u}{\partial \Omega_{\sgn{\tau}}} + \LtOn{u} \\
&\quad\lesssim \q|\tau\w|^{-1} \LtOn{\sLt u} + \q|\tau\w|^{-\frac{1}{2}} \Ltnp{\gamma \partial_{\nu} u}{\partial \Omega_{-\sgn{\tau}}},
\end{split}
\end{equation}
and the same inequality holds with $\sLt$ replaced by $\sLtb$.

Define
$$\sD_{\pm} = \q\{ v\in C^2\q( \overline{\Omega} \w) :  v|_{\partial\Omega}=0, \partial_\nu v|_{\partial\Omega_{\pm}}=0 \w\}.$$
The goal of this section is to prove the following theorem.
\begin{thm}\label{ThmGreensThm}
For any $\tau\ne 0$, there exist operators
$$\Gt, \Gtb:\LtO\rightarrow \LtO$$
such that:
\begin{enumerate}
\item[\rm (i)] $\sLt \Gt = I = \sLtb \Gtb$, i.e., $\lap z^{\tau} \Gt z^{-\tau}=I=\lap \zb^{\tau} \Gtb \zb^{-\tau}$.
\item[\rm (ii)] $\LtOOpn{\Gt}, \LtOOpn{\Gtb} = O\q(\frac{1}{\q|\tau\w|} \w)$ for $\q|\tau \w|>>0$.
\item[\rm (iii)] $\Gt: \LtO\rightarrow z^{-\tau} \HDO$ and for all $u\in \LtO$,
$\tr{\Gt u}$ is supported in $\partial\Omega_{\sgn{\tau}}$.  Similarly,
$\tr{\Gtb u}$ is supported in $\partial \Omega_{\sgn{\tau}}$.
\item[\rm (iv)] $\Gt^{*}=\Gmtb$.
\item[\rm (v)] If $v\in \sD_{-\sgn{\tau}}$, then $\Gt\sLt v=v$, with a similar result for $\Gtb$.
\end{enumerate}
\end{thm}

Let $1-\pt$ be the orthogonal projection onto the closure
in $\LtO$ of $\sLmtb \sD_{\sgn{\tau}}$, and $1-\ptb$ the projection
onto the closure in $\LtO$ of $\sLmt \sD_{\sgn{\tau}}$.  

\begin{lemma}\label{LemmaProj}
$\pt$ is the orthogonal projection onto 
$$\q\{u\in \LtO :  \sLt u =0 \text{ and } \tr{u} \text{ is supported in } \partial\Omega_{\sgn{\tau}}\w\},$$
with a similar result for $\ptb$.
\end{lemma}                                                                     
\begin{proof}                                                                   
Indeed, we will show that $u$ is orthogonal to $\sLmtb\sD_{\sgn{\tau}}$         
if and only if $u$ is as in the statement of the lemma.  Suppose $u$ is         
orthogonal to $\sLmtb\sD_{\sgn{\tau}}$.  Then, in particular,                   
for all $v\in C_0^\infty\q( \Omega\w)$, we have:                                
\begin{equation*}                                                               
\LtOip{u}{\sLmtb v}=0                                                           
\end{equation*}                                                                 
and thus $\sLt u =0$.  Next, allowing $v\in \sD_{\sgn{\tau}}$ to be arbitrary   we see that:                                                                    \begin{equation*}                                                               
0=\LtOip{\sLmtb v}{u} = \int_{\partial \Omega_{-\sgn{\tau}}} \q(\partial_{\nu} v \w) \overline{\tr{u}}                                                          
\end{equation*}
and it follows that $\tr{u}$ is supported in $\partial\Omega_{\sgn{\tau}}$.                                                                                     

The converse follows by the same integration by parts, and is left to the reader.
\end{proof}

The following
lemma yields a unique solution of $\sLt u=f$ which
vanishes on $\partial\Omega_{\sgn{\tau}}$ and is in the range of
$1-\pt$.  The proof is a simple modification of arguments in
\cite{LM}, \cite{BukhgeimUhlmannRecoveringAPotentialFromPartialCauchyData},
\cite{KenigSjostrandUhlmannTheCalderonProblemWithPartialData}.
\begin{lemma}\label{LemmaConstSt}
Given $f\in \LtO$ and $\tau\ne 0$ there exists a unique $u\in \LtO$ such that
\begin{enumerate}
\item\label{NewLemmaItem1} $\sLt u =f,$
\item\label{NewLemmaItem2} $\tr{u}$ is supported in $\partial\Omega_{\sgn{\tau}}$,
\item\label{NewLemmaItem3} $\pt u =0$.
\end{enumerate}
Moreover, this $u$ satisfies $\LtOn{u} \lesssim \q|\tau \w|^{-1} \LtOn{f}$.
\end{lemma}
\begin{proof}
We begin by showing uniqueness.  If $\widetilde{u}$ is another solution
of $\sLt \widetilde{u}=f$ satisfying conditions \ref{NewLemmaItem2} and 
\ref{NewLemmaItem3} above, then $\sLt\q( \widetilde{u}-u\w)=0$ and
$\tr{\widetilde{u}-u}$ is supported in $\partial\Omega_{\sgn{\tau}}$.
Thus, in view of Lemma \ref{LemmaProj}, $\widetilde{u}-u=\pt\q(\widetilde{u}-u\w)$.  However, condition \ref{NewLemmaItem3} above shows $\pt\q( \widetilde{u}-u\w)=0$, and it follows that $\widetilde{u}=u$.

To show existence, define a linear function $l$ initially on $\sLmtb \sD_{\sgn{\tau}}$ by:\footnote{Here 
we are using $\q<\cdot,\cdot\w>_{L^2}$ to denote the {\it sesquilinear} pairing between two $L^2$ functions.  This is in contrast to our notation $\q<\cdot,\cdot\w>$ without the subscript, which denotes the {\it bilinear} pairing of a distribution and a test function.}
$$l\q( \sLmtb v\w) = \LtOip{v}{f}$$
We have:
\begin{equation}\label{EqnEstConsSt}
\q| l\q( \sLmtb v  \w) \w| \leq \LtOn{v}\LtOn{f} \lesssim \q|\tau \w|^{-1} \LtOn{\sLmtb v} \LtOn{f},
\end{equation}
where we have applied the Carleman estimate \eqref{EqnCarleman}.
The functional $l$ extends by continuity to the closure of
$\sLmtb \sD_{\sgn{\tau}}$.  Define $l\equiv 0$ on the orthogonal
complement in $\LtO$ of $\sLmtb \sD_{\sgn{\tau}}$.
There exists a unique $u\in \LtO$ such that:
\begin{equation*}
\LtOip{\sLmtb v}{u} = l\q( \sLmtb v\w) = \LtOip{v}{f}, \quad \q( 1-\pt \w)u=u.
\end{equation*}
Moreover \eqref{EqnEstConsSt} shows that $\LtOn{u}\lesssim \q| \tau\w|^{-1} \LtOn{f}$.
Taking $v\in C_0^{\infty}\q( \Omega\w)$ in the above equation shows that
$\sLt u=f$.  Now
letting $v\in \sD_{\sgn{\tau}}$ be arbitrary, we see via Green's
formula:
\begin{equation*}
\LtOip{\sLmtb v}{u} = \int_{\partial\Omega_{-\sgn{\tau}}} \q(\partial_\nu v\w) \overline{\tr{u}} + \LtOip{v}{f}
\end{equation*}
and therefore, $\int_{\partial\Omega_{-\sgn{\tau}}} \q(\partial_\nu v\w) \overline{\tr{u}}=0$.  Thus, as $v\in \sD_{\sgn{\tau}}$ was arbitrary, 
$\tr{u}$ must be supported in $\partial\Omega_{\sgn{\tau}}$.
\end{proof}

Define $\Ht:\LtO\rightarrow \LtO$ by $\Ht f=u$, where $u$ and $f$ are as
in Lemma \ref{LemmaConstSt}.  In a similar manner,
we construct $\Htb$.  We then have:
\begin{enumerate}
\item $\sLt \Ht = I = \sLtb \Htb$
\item $\q( 1- \pt\w) \Ht =\Ht$, $\q( 1-\ptb\w) \Htb =\Htb$
\item $\LtOOpn{\Ht}=O\q( \q| \tau\w|^{-1}\w)$, $\LtOOpn{\Htb}=O\q( \q| \tau\w|^{-1}\w)$
\item For all $u\in \LtO$, $\tr{\Ht u}$ and $\tr{\Htb u}$ are supported in $\partial \Omega_{\sgn{\tau}}$.
\end{enumerate}
Moreover, $\Ht$ is characterized by the fact that $\q( 1-\pt\w) \Ht=\Ht$ and
\begin{equation*}
\LtOip{\Ht f}{\sLmtb v} = \LtOip{f}{v},\quad \forall v\in \sD_{\sgn{\tau}}.
\end{equation*}
Thus, the operators $\Ht, \Htb$ satisfy (i)-(iii) of Theorem \ref{ThmGreensThm}.
We need to suitably modify $\Ht, \Htb$ to obtain the crucial property (iv).

As a preliminary step, we define
\begin{equation*}
\Tt = \Ht\q( 1-\pmtb\w), \quad \Ttb = \Htb\q( 1-\pmt\w)
\end{equation*}

\begin{lemma}\label{LemmaTtAdj}
$$\Tt^{*} = \Tmtb$$
\end{lemma}
\begin{proof}
Since $\Tt^{*}\pt = 0 = \Tmtb\pt$, it suffices to show that
$$\Tt^{*} \sLmtb v = \Hmtb \sLmtb v, \quad \forall v\in \sD_{\sgn{\tau}}$$
Moreover, since $\q( 1-\pmtb\w) \Tt^{*} = \Tt^{*}$ and $\q( 1-\pmtb\w) \Hmtb = \Hmtb$, it suffices to show:
\begin{equation*}
\LtOip{\sLt w}{\Tt^{*} \sLmtb v} = \LtOip{\sLt w}{\Hmtb \sLmtb v}, \quad \forall w\in \sD_{-\sgn{\tau}}, v\in \sD_{\sgn{\tau}}
\end{equation*}
By the definition of $\Hmtb$ (since $w\in \sD_{-\sgn{\tau}}$), we have:
\begin{equation*}
\LtOip{\sLt w}{\Hmtb \sLmtb v}= \LtOip{w}{\sLmtb v}.
\end{equation*}
We also have:
\begin{equation*}
\begin{split}
\LtOip{\sLt w}{\Tt^{*} \sLmtb v} &= \LtOip{\sLt w}{\q( 1- \pmtb\w) \Ht^{*} \sLmtb v}\\
&= \LtOip{\Ht \q( 1- \pmtb \w) \sLt w}{\sLmtb v}\\
&= \LtOip{\Ht \sLt w}{\sLmtb v}\\
&= \LtOip{\sLt w}{v}\\
&= \LtOip{w}{\sLmtb v}
\end{split}
\end{equation*}
where in the second to last line we used the definition of $\Ht$ and
in the last line we integrated by parts and used the fact that
$w\in \sD_{-\sgn{\tau}}$, $v\in \sD_{\sgn{\tau}}$.  This completes the proof
of the lemma.
\end{proof}

\begin{proof}[Proof of Theorem \ref{ThmGreensThm}]
Define
\begin{equation*}
\Gt = \Ht + \pt \Hmtb^{*}, \quad \Gtb = \Htb + \ptb \Hmt^{*}.
\end{equation*}
It follows from the construction of $\Ht$ and Lemma \ref{LemmaProj}
that $\sLt \Gt =I$, $\tr{\Gt u}$ is supported in $\partial \Omega_{\sgn{\tau}}$
for all
$u\in \LtO$, $\LtOOpn{\Gt}=O\q( \q|\tau\w|^{-1}\w)$, and similar
results hold for $\Gtb$.  To show that $\Gt^{*}=\Gmtb$, we use Lemma \ref{LemmaTtAdj}:
\begin{equation*}
\begin{split}
\Gt^{*} &= \Ht^{*} + \Hmtb\pt\\
&= \q( \Ht \pmtb + \Tt\w)^{*} + \Hmtb - \Tmtb\\
&= \pmtb \Ht^{*} + \Tmtb + \Hmtb - \Tmtb \\
&= \Hmtb + \pmtb \Ht^{*}\\
&= \Gmtb.
\end{split}
\end{equation*}

To verify (v), consider, for $h\in \LtO$, $v\in \sD_{-\sgn{\tau}}$,
\begin{equation*}
\begin{split}
\LtOip{h}{\Gt \sLt v} &= \LtOip{\Gmtb h}{\sLt v}\\
&=\LtOip{\sLmtb \Gmtb h}{v}\\
&=\LtOip{h}{v},
\end{split}
\end{equation*}
where, in the second to last line, we have integrated by parts and
used that $\tr{\Gmtb h}$ is supported in $\partial\Omega_{-\sgn{\tau}}$.
\end{proof}

\section{Special Solutions when $q=0$}\label{SectionSpecialSolutions}
In this section we consider only $\tau>0$.  All of the results in this
section hold for $\tau<0$, provided one reverses
the roles of $\partial \Omega_{+}$ and $\partial\Omega_{-}$ everywhere.
Recall that $\Bt$ is a neighborhood of $\partial\Omega_{+}$.  
The goal of this section is to construct a family of harmonic functions $u_\tau$
in $\Omega$ which vanish on $\partial\Omega\setminus \partial \Bt$ and
have specified asymptotics for large $\tau$.  More precisely, given
any $g\in C^\infty\q( S^{n-2}\w)$, we construct $\mu_\tau=z^{-\tau} u_\tau\in L^2\q( \Omega\w)$ such that:
\begin{enumerate}
\item $\sLt \mut =0$, ie $\lap \ut=0$.
\item The support of $\tr{\mut}$ is in $\Bt$.
\item $\mut\rightarrow \q( z-\zb \w)^{-\frac{n-2}{2}} g\q( \theta\w)$, in $\LtO$ norm as $\tau\rightarrow \infty$.
\end{enumerate}
To do this, we use an extension of Lemma \ref{LemmaConstSt} for
the case where the solution is prescribed
(not necessarily zero) on $\partial\Omega_{-}$.

Define $\sD=\q\{ f\in C^2\q( \Omega\w) : f|_{\partial\Omega}=0 \w\}$,
and define
\begin{equation*}
\Mt=\q\{ \q( \sLmtb f, \partial_{\nu} f|_{\partial\Omega_{+}} \w): f\in \sD \w\}.
\end{equation*}
We think of $\Mt$ as a (non-closed) subspace of $\LtO\times L^2\q( \tau \gamma^2 dS, \partial\Omega_{+} \w)$, where $dS$ denotes the surface
measure on $\partial\Omega$.  Let $\sMt$ denote the orthogonal projection
onto the closure of $\Mt$ in this Hilbert space.  We then have the
following result, which is essentially Proposition 7.1 of \cite{KenigSjostrandUhlmannTheCalderonProblemWithPartialData}.
\begin{prop}\label{PropSpecifyBoundary}
Given 
$$v\in \LtO, v_{-}\in L^2\q(\frac{1}{\gamma^2} dS,\partial\Omega_{-}\w),$$
there exists a unique $u\in \LtO$ such that
\begin{enumerate}
\item $\sLt u = v,$
\item $\tr{u}|_{\partial\Omega_{-}} = v_{-},$
\item $\sMt \q( u, -\tr{u}|_{\partial\Omega_{+}}\w) = \q( u, -\tr{u}|_{\partial\Omega_{+}}\w).$
\end{enumerate}
This $u$ satisfies
$$\LtOn{u} \lesssim \frac{1}{\tau} \LtOn{v} + \tau^{-\frac{1}{2}} \LtdOmn{\frac{1}{\gamma} v_{-}}.$$
\end{prop}
\begin{proof}
Define a linear function $l$ on $\Mt$ by:
$$l\q( \sLmtb f, \partial_{\nu} f|_{\partial\Omega_{+}} \w) = \LtOip{f}{v} + \Ltipp{\partial_{\nu}f}{v_{-}}{\LtdOm}.$$
Let $l\equiv 0$ on the orthogonal complement of $\Mt$.
Note that (using \eqref{EqnCarleman}):
\begin{equation*}
\begin{split}
&\q| l \q( \sLmtb f, \partial_{\nu} f|_{\partial\Omega_{-}}\w)\w| \leq \LtOn{v} \LtOn{f} + \LtdOmn{\gamma \partial_{\nu} f} \q\| v_{-} \w\|_{L^2\q(\frac{1}{\gamma^2} dS,\partial\Omega_{-}\w)}\\
&\quad \lesssim  \q( \frac{1}{\tau} \LtOn{v} + \tau^{-\frac{1}{2}} \LtdOmn{\frac{1}{\gamma} v_{-}} \w) 
\q( \LtOn{\sLmtb f} + \tau^{\frac{1}{2}} \LtdOpn{\gamma \partial_{\nu}f}\w)\\
\end{split}
\end{equation*}
However, the norm applied to $\q( \sLmtb f,\partial_{\nu} f|_{\partial\Omega_{+}}\w)$ on the RHS is precisely the norm in $\LtO\times L^2\q( \tau \gamma^2 dS, \partial\Omega_{+} \w)$.  Thus, there exists 
$$\q( u,u_{+}\w)\in \LtO\times L^2\q( \tau^{-1} \gamma^{-2} dS, \partial\Omega_{+} \w)$$
with $$\LtOn{u} \lesssim \frac{1}{\tau} \LtOn{v} + \tau^{-\frac{1}{2}} \LtdOmn{\frac{1}{\gamma} v_{-}}$$
and satisfying:
\begin{equation*}
\LtOip{\sLmtb f}{u} + \Ltipp{\partial_{\nu}f}{u_{+}}{\LtdOp} = \LtOip{f}{v} + \LtdOmip{\partial_{\nu} f}{v_{-}}
\end{equation*}
Taking $f\in C_0^{\infty}\q( \Omega\w)$ shows that $\sLt u = v$.  Now allowing
$f\in \sD$ to be arbitrary shows that $\tr{u}|_{\partial\Omega_{-}}=v_{-}$
and $\tr{u}|_{\partial\Omega_{+}}=-u_{+}$.  It is clear that this $u$
is unique under the conditions of the proposition.
\end{proof}
Let $\Rt \q( v, v_{-}\w) = u$, where $u, v, v_{-}$ are as in Proposition \ref{PropSpecifyBoundary}.
By reversing the roles of $\partial\Omega_{+}$ and $\partial\Omega_{-}$
we also get $\Rmt \q( v, v_{+}\w)$ where $v_{+}\in L^2\q(\frac{1}{\gamma^2} dS,\partial\Omega_{+}\w)$.

We now turn to the construction of the solutions promised at the beginning
of this section.  It is easy to see, using \eqref{EqnLapEqn}, that:
\begin{equation*}
\begin{split}
\sLt &= \lap +\tau\q( \frac{4}{z} \frac{\partial}{\partial \zb} - \frac{2\q( n-2\w)}{\q( z-\zb \w) z} \w)\\
&= \lap + \frac{\tau}{z} L
\end{split}
\end{equation*}
where
$$L=4\frac{\partial}{\partial \zb} - \frac{2\q( n-2\w)}{\q( z-\zb\w)}$$
Note that:
$$L\q( \q(z-\zb\w)^{-\frac{n-2}{2}}g\q( \theta\w)  \w)=0$$
and so:
$$\sLt \q(z-\zb\w)^{-\frac{n-2}{2}}g\q( \theta\w) = \lap \q(z-\zb\w)^{-\frac{n-2}{2}}g\q( \theta\w) $$
Hence, if we let $\chi_{+}\in C_0^\infty\q( \q\{ x\in \partial\Omega : x\cdot \nu<0 \w\}\w)$, with $\chi_{+}=1$ on an open set containing $\partial\Omega\setminus \Bt$, 
we see that (if $h_{+} = \q(z-\zb\w)^{-\frac{n-2}{2}}g\q( \theta\w)$):
$$\mut = h_{+} - \Rt\q( \lap h_{+}, \q( \chi_{+}\w) \tr{h_{+}} \w)$$
has the desired properties.

In addition, if we take $\chi_{-}\in C_0^\infty\q( \q\{ x\in \partial\Omega : x\cdot \nu>0 \w\}\w)$, with $\chi_{-}=1$ on a neighborhood of $\partial\Omega \setminus \Ft$, and
if we define $h_{-}= \q(z-\zb\w)^{-\frac{n-2}{2}}$,
$$\numt = h_{-} - \Rmt\q( \lap h_{-}, \q( \chi_{-}\w) \tr{h_{-}}\w)$$
satisfies $\sLmt {\numt} =0$, the support of $\tr{\numt}$ is in $\Ft$,
and $\numt\rightarrow \q( z- \zb\w)^{-\frac{n-2}{2}}$ in 
$\LtO$ norm as $\tau \rightarrow \infty$.  For the
rest of the paper, we fix this choice of harmonic functions $u_{\tau}=z^\tau \mut$ and $\vmt = z^{-\tau}\numt$.

\section{Special Solutions for General $q$}\label{SectionSpecialSolutionsq}
In this section we construct our family of solutions $\wt=z^\tau \ot$
of $\q( -\lap +q\w) \wt =0$ in $\Omega$, which vanish on $\partial\Omega\setminus \Bt$ and have specified asymptotics for large $\tau$.

We also define a corresponding nonlinear transform of $q$.  We take
$\ut$ and $\vmt$ as constructed in Section \ref{SectionSpecialSolutions}.
Recall that $\ut$ was defined in terms of a fixed $g\in C^\infty\q( S^{n-2}\w)$.

\begin{prop}\label{PropConstot}
Let $\Omega$, $q$, $x_0$, $F\q( x_0\w)$, $B\q( x_0\w)$, $\Ft$, and $\Bt$
be as in the hypothesis of Theorem \ref{ThmMainThm}.  For $\tau>>0$,
$g\in C^\infty\q( S^{n-2}\w)$ and $\mut$ as above, there exists a
unique solution $\ot=z^{-\tau} w_\tau$ of the integral equation
\begin{equation}\label{EqnIntEqn}
\ot = \mut+\Gt q \ot.
\end{equation}
This solution satisfies
\begin{enumerate}
\item[{\rm (i)}] $\q( -\sLt +q \w) \ot =0$,
\item[{\rm (ii)}] $\tr{\ot}$ is supported in $\Bt$,
\item[{\rm (iii)}] $\ot\rightarrow \q( z-\zb\w)^{\frac{n-2}{2}}g\q( \theta\w)$ in $\LtO$ norm as $\tau\rightarrow \infty$.
\end{enumerate}
\end{prop}
\begin{proof}
Since $q\in L^\infty\q( \Omega\w)$, unique solvability of \eqref{EqnIntEqn}
for $\tau$ sufficiently large follows from the bound
$\LtOOpn{\Gt}=O\q( \tau^{-1}\w)$ (Theorem \ref{ThmGreensThm} (ii)).
Property (i) follows from $\sLt\mut=0$ and $\sLt \Gt=I$.  Property
(ii) is a consequence of the support properties of
$\tr{\mut}$ and $\trnp\circ \Gt$ (Theorem \ref{ThmGreensThm} (iii)).
Finally, (iii) follows from the corresponding asymptotics of $\mut$
and the above bound on $\Gt$.
\end{proof}

By analogy with the approach to Calder\'on's problem
in \cite{N1}, \cite{N2}, we define the following nonlinear transform
of $q$
\begin{equation}\label{EqnDefnt}
t\q( \tau, g\w) = t\q( \tau, g, x_0, \alpha_0\w) =\int_\Omega \vmt q \wt = \int_\Omega \numt q \ot
\end{equation}
with $\ot$ as constructed above and $\numt$ the solution for
homogeneous background defined in Section \ref{SectionSpecialSolutions}.
Then
\begin{equation}\begin{split}\label{EqnLimitTt}
\lim_{\tau\rightarrow \infty} t\q( \tau,g\w) &= \int_\Omega \q( z-\zb\w)^{-\frac{n-2}{2}} q \q( z-\zb\w)^{-\frac{n-2}{2}} g
\\
&= \int_\Omega q\q( x_1,r,\theta\w) g\q( \theta\w) \q( 2ir\w)^{-\q( n-2\w)} r^{n-2}\: dr\: d\theta\: dx_1
\\
&= \q( 2i\w)^{-\q( n-2\w)}\int_\Omega q\q( x_1, r, \theta\w) g\q( \theta\w) \: dr\: d\theta\: dx_1.
\end{split}
\end{equation}
Thus, to prove Theorem \ref{ThmMainThm} it suffices to reconstruct
$t\q( \tau, g\w)$ from the given partial knowledge of the
Dirichlet-to-Neumann map.
\begin{thm}\label{ThmReconTt}
Given partial knowledge of $\Lambda_q$ as in the hypothesis of Theorem
\ref{ThmMainThm}, one can reconstruct $t\q( \tau, g\w)$ for any
$g\in C^\infty\q( S^{n-2}\w)$ and $\tau$ sufficiently large.
\end{thm}
\begin{proof}\renewcommand{\qedsymbol}{}
Using \eqref{EqnAlles} we can express $t\q( \tau, g\w)$ in terms of boundary
data as
\begin{equation}\label{EqnTtEquals}
t\q( \tau, g\w) = \int_\Omega \vmt q \wt = \int_{\partial \Omega} \tr{\vmt} \q( \Lambda_q -\Lambda_0\w) \tr{\wt}.
\end{equation}
Recall that $\vmt$ was defined independently of $q$.  Since $\tr{\vmt}$
is supported in $\Ft$ and $\tr{\wt}$ is supported in $\Bt$, the above
formula only involves the given partial knowledge of $\Lambda_q$.  The
proof will be completed in the following section, where we will show that
$\tr{\wt}$ can be reconstructed from the given data.
\end{proof}

\section{Single Layer Operators and the Boundary Integral Equation}\label{SectionIntegralEqns}
In this section, we define the single layer operators $\St$ corresponding
to the Green's operators constructed in Section \ref{SectionGreensOp}
and show that $\tr{\wt}$ can be reconstructed as
the unique solution of the integral equation in $\sH$:
\begin{equation}\label{Eqn105}
\tr{\wt} = \tr{\ut} + \St \q( \Lambda_q-\Lambda_0\w) \tr{\wt}
\end{equation}
(Compare with (0.17) in \cite{N2}.)

Consider the map
$$\q(\trnp\circ \Gt \w)^{*} : \zb^\tau \sH^{*}\rightarrow \LtO$$
defined for $h\in \zb^{\tau} \sH^{*}$, $f\in \LtO$ by:
$$\int_\Omega \q(\overline{\q(\trnp\circ \Gt\w)^{*} h}\w) f = \int_{\partial\Omega}\overline{h} \q(\trnp\circ \Gt\w)f$$
\begin{lemma}\label{Lemma61}
For all $h\in \zb^\tau \sH^{*}$,
\begin{enumerate}
\item[{\rm (i)}] $\sLmtb\q(\trnp\circ \Gt\w)^{*} h=0$
\item[{\rm (ii)}] $\tr{\q(\trnp\circ \Gt\w)^{*}h}$ is supported in $\partial \Omega_{-}$, if $\tau>0$.
\end{enumerate}
Moreover if $\Btt$ is a neighborhood of $\partial\Omega_{+}$ such that
$\overline{\Btt}\subset \Bt$, and if $h$ is supported in $\partial\Omega\setminus \Btt$, we have that $\q(\trnp\circ \Gt\w)^{*} h=0$.
\end{lemma}
\begin{proof}
If $h$ is supported in $\partial\Omega\setminus \Btt$,
$$\int_\Omega \q(\overline{\q(\trnp\circ \Gt\w)^{*} h}\w) f = \int_{\partial\Omega_{+}}\overline{h} \q(\trnp\circ \Gt\w) f=0$$
for all $f\in \LtO$, since $\q( \trnp\circ \Gt\w) f$ is supported in $\partial \Omega_{+}$.  Thus, $\q( \trnp\circ \Gt\w)^{*}h=0$.

Now consider, for $f\in \sD_{-}$ and all $h$,
\begin{equation}\label{Eqn106}
\begin{split}
\int_\Omega \q(\overline{\q(\trnp \circ \Gt\w)^{*} h}\w) \sLt f &= \int_{\partial\Omega} \overline{h} \q( \trnp\circ \Gt\w) \sLt f\\
&= \int_{\partial\Omega} \overline{h} \tr{f}\\
&=0
\end{split}
\end{equation}
where in the second to last line we have used Theorem \ref{ThmGreensThm} (v),
and in the last line we have used that $\tr{f}=0$.  The rest of the lemma now
follows from \eqref{Eqn106} and integration by parts (similar to that in Lemma \ref{LemmaConstSt}).
\end{proof}

For any $\tau\ne 0$, define the operator $\St$ on $\sH^{*}$ as
\begin{equation}\label{Eqn107}
\St h = z^\tau \q( \trnp\circ\q(\trnp\circ \Gt\w)^{*} \w)^{*}\q(hz^{-\tau}\w)
\end{equation}
To streamline the notation we assume $\tau>0$ in some of the results below.

\begin{prop}\label{Prop62}
\begin{enumerate}
\item[{\rm (i)}] $S_\tau$ is a bounded operator $\sH^{*}\rightarrow \sH$.
\item[{\rm (ii)}] If $\tau>0$, $\St h$ only depends on $h\big|_{\Ft}$.
\item[{\rm (iii)}] If $\tau>0$, $\St h$ is supported in $\Bt$.
\end{enumerate}
\end{prop}
\begin{proof}
It follows from Lemma \ref{Lemma61} (i) that
$$\q(\trnp\circ \Gt\w)^{*}: \zb^\tau \sH^{*} \rightarrow \zb^\tau \HDO,$$
hence $\trnp\circ\q(\trnp \circ \Gt\w)^{*}: \zb^\tau \sH^{*} \rightarrow \zb^{\tau} \sH$.  Thus, $\q(\trnp\circ\q(\trnp\circ \Gt\w)^{*}\w)^{*}$ is a bounded
operator $z^{-\tau}\sH^{*} \rightarrow z^{-\tau} \sH$ and (i) follows.

\noindent (ii)  By Lemma \ref{Lemma61} (ii), if $\htt\in \sH^{*}$ then
$\trnp\circ\q( \trnp\circ\Gt\w)^{*} \zb^{\tau} \htt$ is supported
in $\partial \Omega_{-}$, so that for $h\in \sH^{*}$ one may compute
$z^\tau \q(\trnp\circ\q(\trnp\circ \Gt\w)^{*}\w)^{*} z^{-\tau} h$
using only knowledge of $h|_{\Ft}$.

\noindent (iii) Since $\trnp\circ\q(\trnp\circ \Gt\w)^{*}\zb^{\tau} \htt=0$
for all $\htt\in \sH^{*}$ supported in $\partial\Omega\setminus \Btt$,
we see that $z^{\tau} \q(\trnp\circ\q(\trnp\circ \Gt\w)^{*}\w)^{*}z^{-\tau}h$
is supported in $\Bt$ for any $h\in \sH^{*}$.
\end{proof}

It follows from the above and Proposition \ref{PropExtendDTN} that the operator
$\St \q( \Lambda_q - \Lambda_0\w)$ is bounded $\sH\rightarrow \sH$,
and that for $\tau>0$ and any $h\in \sH$,
$\St\q( \Lambda_q-\Lambda_0\w) h$ is supported in $\Bt$.  Moreover,
if $\tau>0$ and $h$ is supported in $\Bt$, then $\St\q( \Lambda_q-\Lambda_0\w)h$
can be computed using only the given partial knowledge of $\Lambda_q$.

To prove the solvability of the integral equation \eqref{Eqn105},
we'll use the following factorization identity (compare with 
(7.4) of \cite{N2}).
\begin{prop}\label{Prop63}
$$\St \q( \Lambda_q-\Lambda_0\w) = \trnp\circ z^{\tau} \Gt z^{-\tau} q \sPqnp.$$
\end{prop}
\begin{proof}
Consider, for any $f\in \sH^{*}$, $h\in \sH$,
$$\q< \overline{f},\trnp\circ z^{\tau} \Gt z^{-\tau} q \sPqnp h\w> = \int_{\Omega} \q(\overline{\zb^{-\tau} \q(\trnp\circ \Gt\w)^{*} \zb^{\tau} f}\w) q\sPqnp h.$$

By Lemma \ref{Lemma61}, $\zb^{-\tau} \q(\trnp\circ \Gt\w)^{*} \zb^\tau f\in \bz$ and,
using \eqref{EqnAlles} we find that the right hand side of the above equality
equals
\begin{equation*}
\begin{split}
\q<\overline{\zb^{-\tau}\trnp\circ\q(\trnp\circ \Gt\w)^{*} \zb^{\tau}f},\q(\Lambda_q-\Lambda_0\w)h\w> &= \q<\overline{f}, z^\tau \q(\trnp\circ\q(\trnp\circ \Gt\w)^{*}\w)^{*}z^{-\tau} \q(\Lambda_q-\Lambda_0\w)h\w>
\\&= \q<\overline{f}, \St \q(\Lambda_q-\Lambda_0\w)h\w>.
\end{split}
\end{equation*}
\end{proof}

The next three results below will yield the solvability of \eqref{Eqn105}.

\begin{prop}\label{Prop64}
For $f,h\in \sH$,
$$\q[I-\St\q( \Lambda_q-\Lambda_0\w)\w]h=f$$
if and only if
$$\q(I-z^{\tau} \Gt z^{-\tau} q\w)\sPqnp h = \sPznp f.$$
\end{prop}

\begin{cor}\label{Cor65}
$$I-\St\q( \Lambda_q-\Lambda_0\w): \sH\rightarrow \sH$$
 is an
isomorphism if and only if
$$I-z^{\tau} \Gt z^{-\tau} q: \bq\rightarrow \bz$$
is an isomorphism.
\end{cor}

\begin{prop}\label{Prop66}
For $\q|\tau\w|>>0$, $$I-z^{\tau} \Gt z^{-\tau} q:\bq\rightarrow \bz$$
is an isomorphism.
\end{prop}

\begin{proof}[Proof of Proposition \ref{Prop64}]
Suppose $\q[I-\St \q(\Lambda_q-\Lambda_0\w)\w]h=f$.  We wish to show that
$$\q(I-z^\tau \Gt z^{-\tau} q\w)\sPqnp h = \sPznp f.$$
Since
$$\lap\q(I-z^\tau \Gt z^{-\tau} q\w)\sPqnp h = q\sPqnp h- q\sPqnp h =0,$$
it suffices to show that
$$\tr{\q(I-z^\tau \Gt z^{-\tau} q\w)\sPqnp h} = f$$
ie that
$$h-\tr{z^\tau \Gt z^{-\tau} q\sPqnp h} =f.$$
However, Proposition \ref{Prop63} shows that the left side equals
$h-\St\q(\Lambda_q-\Lambda_0\w)h$, which we are assuming equal to $f$.

For the converse, suppose that
$$\q(I-z^\tau \Gt z^{-\tau}q\w)\sPqnp h = \sPznp f$$
taking the trace of both sides yields, using Proposition \ref{Prop63}
$$h-\St\q(\Lambda_q-\Lambda_0\w)h=f,$$
as claimed.
\end{proof}

\begin{proof}[Proof of Proposition \ref{Prop66}]
Since $q\in L^\infty\q( \Omega\w)$ and since $\LtOOpn{\Gt}=O\q( \frac{1}{\q|\tau\w|}\w)$ (Theorem \ref{ThmGreensThm} (ii)),
$I-\Gt q:\LtO\rightarrow \LtO$ is an isomorphism for $\q|\tau\w|$
sufficiently large.
For such $\tau$, the operator $I-z^\tau \Gt z^{-\tau}q$ is then
an isomorphism $\LtO\rightarrow \LtO$ with
inverse $z^{\tau}\q(I-\Gt q\w)^{-1} z^{-\tau}$.  If $u\in \bz$,
we claim 
$$w=z^{\tau} \q(I-\Gt q\w)^{-1} z^{-\tau} u\in \bq.$$  
Indeed,
$$w-z^{\tau} \Gt q z^{-\tau} w =u$$
and $\lap z^\tau \Gt z^{-\tau} =I$, hence $\lap w- qw=0$.
\end{proof}

\begin{proof}[Proof of Theorem \ref{ThmReconTt} (completed)]
Let $\wt$ be defined as in Proposition \ref{PropConstot}.
Then, in view of Proposition \ref{Prop64}, $\tr{\wt}$
satisfies \eqref{Eqn105}.  Corollary \ref{Cor65} and
Proposition \ref{Prop66} show that \eqref{Eqn105} is
uniquely solvable for large $\tau$.  Substituting the solution
$\tr{w_\tau}$ in formula \eqref{EqnTtEquals} yields
$t\q( \tau, g\w)$.
\end{proof}

\begin{proof}[Proof of Theorem \ref{ThmMainThm}]
We solve the boundary integral equation \eqref{Eqn105} for
$\tau>>0$, as indicated above.  We then calculate
$t\q( \tau, g\w)$ using formula \eqref{EqnTtEquals}.
The large $\tau$ limit \eqref{EqnLimitTt} then gives
the integrals \eqref{EqnIntToRecon}, as claimed.
\end{proof}



\providecommand{\bysame}{\leavevmode\hbox to3em{\hrulefill}\thinspace}
\providecommand{\MR}{\relax\ifhmode\unskip\space\fi MR }
\providecommand{\MRhref}[2]{%
  \href{http://www.ams.org/mathscinet-getitem?mr=#1}{#2}
}
\providecommand{\href}[2]{#2}

\end{document}